\documentclass{amsart}
\usepackage{amssymb}
\usepackage{amsmath}

\newtheorem{thm}{Theorem}[section]

\newtheorem{conj}[thm]{Conjecture}
\newtheorem{question}[thm]{Question}
\newtheorem{problem}[thm]{Problem}

\theoremstyle{definition}
\newtheorem{defn}[thm]{Definition}

\theoremstyle{remark}

\newcommand{\R}{\mathbb{R}}

\newcommand{\calC}{\mathcal C}

\newcommand{\duaug}[2]{{#1}^{+{#2}}}
\newcommand{\link}{\operatorname{Link}}

\newcommand{\calP}{\mathcal{P}}

\usepackage{hyperref}

\input xy
\xyoption{all}

\author{Alessandro Sisto}
	\address{Maxwell Institute and Department of Mathematics, Heriot-Watt University, Edinburgh, UK}

\begin{document}

\title{New tools in hierarchical hyperbolicity: A survey}

\begin{abstract}
 The aim of this short survey is to advertise various tools that have been developed to study hierarchically hyperbolic spaces (HHSs) in recent years, with particular emphasis on those that require little to no knowledge of the HHS machinery to be used.
\end{abstract}

\maketitle

\section{Introduction}

Hierarchically hyperbolic groups were defined in \cite{HHS_I} as a generalisation of hyperbolic groups and a common framework for studying mapping class groups and (special) cube complexes. The definition is an axiomatic version of the machinery developed by Masur and Minsky to study mapping class groups \cite{MM1,MM2}.

As it turns out, there are a lot more hierarchically hyperbolic groups and spaces (HHGs and HHSs) than those motivating the original definition. In fact, examples of a wide variety of different flavours are now known, from ``most'' types of 3-manifold groups to extra-large type Artin groups, see among others \cite{HHS_II,BerlaiRobbio,HagenSusse,hhs_asdim,BHMS,Chesser:stable,Vokes,HRSS_3manifold,ELTAG_HHS,veech,BR:graphs,HV:not_automatic, bongiovanni2024extensions}. Many new results have been obtained using the HHS machinery, see for instance \cite{HHS_I,DHS,hhs_asdim,ABD,HHS:quasiflats,HHS-convexity,ANS:UEG,Petyt_MCG_cubical,HHP:coarse,DMS:stable,induced_qi,Mangioni:combination,barak:eq_not}, and many of these are new even for mapping class groups or other very well-studied examples. Some of these results are mentioned below, just to mention a couple more here, these include the strong Tits alternative \cite{DHS,DHS:oops}, the quasiflat theorem \cite{HHS:quasiflats}, the existence of bicombings \cite{HHP:coarse}, uniform exponential growth for a large subclass \cite{ANS:UEG}, and mapping class groups being quasi-isometric to CAT(0) cube complexes \cite{Petyt_MCG_cubical}.

The definition of hierarchical hyperbolicity is explained and commented on in detail in the survey \cite{HHS_survey}. This survey, instead, focuses on a variety of tools and techniques that have been developed in recent years to study HHSs, and which make the study of HHSs much more accessible to the ``uninitiated'', in various ways. Indeed, the content of some of them is that HHSs also carry other very rich structures, such as injective metrics. The theory of combinatorial HHSs, on the other hand, makes it substantially easier to construct new examples of HHSs, as one no longer needs to check the definition directly. For the purposes of this survey, no particular knowledge of HHSs is required, and we only recall very briefly that the main data of an HHS structure on a metric space $X$ is a collection of uniformly hyperbolic spaces $\{\calC(Y)\}_{Y\in\mathcal S}$, each equipped with a coarsely Lipschitz map $\pi_Y:X\to \calC(Y)$ (with uniform constants across all $Y$). The interested reader is referred to the survey \cite{HHS_survey} for more information, where we note that more ``hierarchical'' and basic tools are discussed already. Notably, the cubulation of hulls from \cite{HHS:quasiflats} is discussed there, and interestingly a few of the results mentioned below (including both re-metrisation results) heavily rely on it and its further refinements and improvements, see \cite{Bowditch:quasiflats,DMS:stable, DZ:genericity, Durham:infinity}. I decided to not discuss this further here as I see this as more intrinsically hierarchical than the tools I want to focus on.

\subsection{Outline}

In Section \ref{sec:combinatorial} we describe combinatorial HHSs, which provide a very useful criterion to construct new HHSs. Indeed, after the development of combinatorial HHSs there has been a substantial increase in activity in finding new HHSs.

In Section \ref{sec:metrics} we describe two ways of re-metrising an HHS. One yields injective metrics, while the other (which works under an additional but very general condition) yields asymptotically CAT(0) metrics. Both are metrics with nonpositive curvature characteristics and have been studied independently of HHSs (especially injective metrics).

In Section \ref{sec:curtains} we discuss various other tools being developed, including Dehn-filling-like quotients, curtains, $\mathbb R$-cubings, and higher rank JSJ decompositions.

\subsection*{Acknowledgements} I would like to thank Jason Behrstock, Mark Hagen, Giorgio Mangioni, Harry Petyt, and Abdul Zalloum for their very useful comments and suggestions.

\section{Combinatorial HHSs}
\label{sec:combinatorial}

In this section we describe the criterion for hierarchical hyperbolicity devised in \cite{BHMS}, which has now been used numerous times to construct new HHSs. It seems likely that more new examples will be constructed using this criterion. The discussion below covers the main idea and formulates the main statement, but more information is available in the ``User's guide'' in \cite{BHMS}, and also in \cite{converse}, where a more explicit description of the candidate combinatorial model for an HHS is given.

First, an informal discussion. Roughly speaking, the criterion says that if $G$ acts ``nicely'' on a hyperbolic simplicial complex $X$ where all links of simplices are also hyperbolic, then $G$ is hierarchically hyperbolic. (By link of a simplex we mean the intersection of the links of its vertices.) Here, ``nicely'' includes the requirements that the action is cocompact and that stabilisers of maximal simplices are finite. With these hypotheses, the group $G$ is quasi-isometric to a graph $W$ whose vertex set is the set of maximal simplices of $X$, and where one adds finitely many orbits of edges to make this set connected. (This is a simple generalization of letting a finitely generated group $G$ act on itself, and adding finitely many orbits of edges to form the Cayley graph.) The idea is then that the hyperbolic spaces for an HHS structure on $G$ are all links of simplices of $X$ (where here we include the link of the empty simplex, which is $X$ itself). This is not accurate, though, because $X$ will have some discrete links, so it is basically never the case that a simplicial complex has all links hyperbolic. To remedy this, one considers ``augmented links'', adding edges suitably. This is similar to the procedure adding edges to the set of maximal simplices to obtain a graph quasi-isometric to the group.

With this discussion in mind, here are the actual definitions. Throughout, let $X$ be a finite-dimensional simplicial complex. First of all, we describe the setup of a pair (simplicial complex, graph whose vertices are maximal simplices), and also we describe additional edges that will be needed later.

\begin{defn}[$X$--graph, $W$--augmented graph]\label{defn:X_graph}
An \emph{$X$--graph} is a graph $W$ whose vertex set is the set of all maximal simplices of $X$.

For a flag simplicial complex $X$ and an $X$--graph $W$, the \emph{$W$--augmented graph} $\duaug{X}{W}$ is the graph defined as follows:
\begin{itemize}
     \item the $0$--skeleton of $\duaug{X}{W}$ is $X^{(0)}$;
     \item if $v,w\in X^{(0)}$ are adjacent in $X$, then they are adjacent in $\duaug{X}{W}$;
     \item if two vertices in $W$ are adjacent, then we consider $\sigma,\rho$, the associated maximal simplices of $X$, and in $\duaug{X}{W}$ we connect each vertex of $\sigma$ to each vertex of $\rho$.
\end{itemize}
\end{defn}

Next, we describe augmented links, that is, links with additional edges. We also need an auxiliary space, denoted $Y_\Delta$ below, which is perhaps the least intuitive feature of the whole setup.

\begin{defn}[Augmented links]\label{defn:complement}
Let $W$ be an $X$--graph.  For each simplex $\Delta$ of $X$, let its \emph{augmented link} $\calC  (\Delta)$ be the induced subgraph of $\duaug{X}{W}$ spanned by $\link(\Delta)^{(0)}$. Also, let $Y_\Delta$ be the subgraph of $\duaug{X}{W}$ induced by the set of vertices
$$X^{(0)}-\bigcup_{\Sigma:\ \link(\Sigma)=\link(\Delta)}\Sigma^{(0)}.$$
\end{defn}

Finally, we describe the axioms of combinatorial HHSs. A more general version of item \eqref{item:chhs_join} below suffices, but we state a simpler one which most often holds in practice. Note that all conditions except for \eqref{item:chhs_delta} are combinatorial conditions on $X$ and $W$, with \eqref{item:chhs_delta} being the only geometric condition.

\begin{defn}[Combinatorial HHS]\label{defn:combinatorial_HHS}
A \emph{combinatorial HHS} $(X,W)$ consists of a flag simplicial complex $X$ and an $X$--graph $W$ satisfying the following conditions:
\begin{enumerate}
    \item There exists $n\in\mathbb N$ such that any chain $\link(\Delta_1)\subsetneq\dots\subsetneq\link(\Delta_i)$, where each $\Delta_j$ is a simplex of $X$, has length at most $n$;
    \item \label{item:chhs_delta} There is a constant $\delta$ so that for each non-maximal simplex $\Delta$, the subgraph $\calC  (\Delta)$ is $\delta$--hyperbolic and $(\delta,\delta)$--quasi-isometrically embedded in $Y_\Delta$;
    \item \label{item:chhs_join} Whenever $\Delta$ and $\Sigma$ are non-maximal simplices, there exist simplices $\Pi,\Pi'$ such that $\link(\Delta)\cap\link(\Sigma)=\link(\Delta\star\Pi)\star \Pi'$, where $\star$ denotes the join;
    \item \label{item:C_0=C} If $v,w$ are distinct non-adjacent vertices of $\link(\Delta)$, for some simplex $\Delta$ of $X$, contained in $W$-adjacent maximal simplices, then they are contained in $W$-adjacent simplices of the form $\Delta\star\Sigma$.
\end{enumerate}
\end{defn}

We are finally ready to state the criterion from \cite{BHMS}.

\begin{thm}\label{thm:hhs_links}
Let $(X,W)$ be a combinatorial HHS. Then $W$ is a hierarchically hyperbolic space.
\end{thm}

Note that \cite{BHMS} actually gives a concrete HHS structure for $W$, and the associated hyperbolic spaces are exactly the augmented links. When there is a group action of $G$ on $(X,W)$ which is compatible in a suitable sense, then one can conclude that $G$ is a hierarchically hyperbolic group.

Very recently, a different version of combinatorial HHS has been used in \cite{Park:factor_system}.

\subsection{Case study: pants graphs} To conclude, we illustrate the criterion in a case where it works incredibly well, that of the pants graph. The pants graph $\calP$ of a fixed surface $\Sigma$ is the graph with vertex set all pants decompositions of $\Sigma$, and where two pants decomposition are connected if the differ by a ``move''. Here, a move consists in changing a single curve with one that intersects it once or twice (depending whether the other curves cut out a one-holed torus or a four-holed sphere). Notice that the vertices of $\calP$ are exactly the maximal simplices of the curve graph $\calC$ of $\Sigma$, so that, in the terminology above, $\calP$ is a $\calC$-graph. Checking the combinatorial conditions is left as an exercise to the reader. Each augmented link is the curve graph $\calC(Y)$ of a subsurface $Y$, with additional edges connecting curves that intersect once or twice, but these do not change the quasi-isometry type. Therefore, augmented links are hyperbolic. Finally, the quasi-isometric embedding part of condition \eqref{item:chhs_delta} can be checked using subsurface projections. As it turns out, for $Y$ the subsurface associated to the simplex $\Delta$, the curves that have well-defined subsurface projection to $\calC(Y)$ correspond exactly to the vertex set of $Y_\Delta$, and it is easy to see that subsurface projections are coarsely Lipschitz (where they are defined). This provides a coarse retraction of $Y_\Delta$ onto $\calC(\Delta)$, therefore showing that the former is quasi-isometrically embedded in the latter. To my knowledge, in all applications of combinatorial HHSs the quasi-isometric embedding condition is checked by constructing a suitable coarse retraction.

\section{Re-metrisation}
\label{sec:metrics}

Hierarchically hyperbolic groups display features of nonpositive curvature, but not all of them are nonpositively curved in the most classical sense, meaning that not all of them are CAT(0). Indeed, most mapping class groups are not CAT(0) \cite{Bridson:MCG}, but also it is not known whether all hyperbolic groups are CAT(0). Despite this, hierarchically hyperbolic groups act properly and coboundedly on spaces that are nonpositively curved in other senses.

\subsection{Injective metrics}  The first result of this type is due to Haettel-Hoda-Petyt \cite{HHP:coarse}, who showed:

\begin{thm}\cite{HHP:coarse}
 All hierarchically hyperbolic groups act properly and coboundedly on some injective space.
\end{thm}

Injective spaces have different characterisations, one being that they are geodesic metric spaces where, given balls that pairwise intersect, they all have a common intersection. Often more useful than the definition is the explicit construction of the \emph{injective hull} of a metric space $X$, which is the smallest injective space containing $X$ (in particular, $E(X)=X$ if $X$ is injective). Injective spaces have various features of nonpositive curvature, including barycenters for finite sets and bicombings with convexity properties.  Both of these were exploited in \cite{HHP:coarse} to prove new results about HHSs, namely that they contain finitely many conjugacy classes of finite subgroups, and that they are semi-hyperbolic, and hence have solvable conjugacy problem.

Injective metric spaces have been brought to the attention of geometric group theorists by Lang, who showed that all hyperbolic groups act geometrically on injective spaces \cite{Lang}, despite the fact that we still do not know whether all hyperbolic groups are CAT(0). They have been an active area of study since, and the theorem above provided further motivation to pursue their study. The reader interested in injective spaces (perhaps to obtain new results about HHSs) is referred to \cite{Haettel:survey,Lang, Zalloum_survey}. As a word of warning, it is not known whether one can improve ``coboundedly'' to ``cocompactly'' in the theorem, and probably one in fact cannot do this. This sometimes causes problems when applying results on injective spaces to HHGs.

\subsection{Asymptotically CAT(0) metrics} There is now another class of spaces for HHGs to act on, even though (at least for the moment) one needs the extra technical assumption of colourability, where there is essentially only one example of non-colourable HHG, constructed specifically to have this property in \cite{noncolorable}. The relevant class of spaces is that of \emph{asymptotically CAT(0) spaces}, which were defined by Kar \cite{Kar:thesis} (at least in the case of geodesic spaces, we need a slight extension here). Roughly, a space is asymptotically CAT(0) if all triangles satisfy the CAT(0) inequality up to an error sublinear in their size. Note that injective spaces have fine nonpositive curvature features, while asymptotically CAT(0) spaces have coarse nonpositive curvature features. Here are all the relevant definitions.

Given 3 points $x,y,z$ in a metric space, we denote $\bar\Delta(x,y,z)$ a comparison triangle in the Euclidean plane $\mathbb E^2$, with vertices $\bar x,\bar y,\bar z$. Given a point $p$ on a $(1,C)$-quasi-geodesic $\gamma$ joining $x,y$, for a fixed $C$, we call a point $\bar p$ on the geodesic $[\bar x,\bar y]$ a \emph{comparison point} if $
|d(p,x)-d(\bar p, \bar x)|\leq C$.
 Also, we will denote $d(p)=\min\{d(x,p),d(p,y)\}$.

\begin{defn}[Sublinear CAT(0)]\label{defn:sub CAT0}
    Given a sublinear and non-decreasing function $\kappa$, we say that a triangle $\Delta$ of $(1,C)$-quasi-geodesics satisfies the \emph{CAT(0) condition up to $\kappa$} if the following holds. Let $x,y,z$ be the vertices of the triangle, and let $p$ and $q$ be points on the triangle. Fixing a comparison triangle and comparison points, we have
$$d(p,q)\leq d(\bar p,\bar q)+\kappa(d(p))+\kappa(d(q)).$$
\end{defn}

The following definition is a small variation on \cite[Definition 6]{Kar:thesis}:

\begin{defn}\label{defn:asymp CAT0}
    Given a constant $C$ and a sublinear function $\kappa$, we say that a metric space $X$ is $(C,\kappa)$-\emph{asymptotically CAT(0)} if the following hold:
    \begin{enumerate}
        \item Every pair of points of $X$ is connected by a $(1,C)$-quasi-geodesic.
        \item Every triangle $\Delta$ of $(1,C)$-quasi-geodesics satisfies the CAT(0) condition up to $\kappa$.
    \end{enumerate}
\end{defn}

One of the main theorems of \cite{FJ} is:

\begin{thm}
 All colourable hierarchically hyperbolic groups act geometrically on some $(C,\kappa)$-asymptotically CAT(0) space.
\end{thm}

The main motivation to prove this theorem is that it is a key step in proving the Farrell-Jones conjecture for a large class of hierarchically hyperbolic groups, which is the main application in \cite{FJ}. This is a bit more restrictive than colourable ones, as we do not quite prove that a geometric action on an asymptotically CAT(0) space suffices for the Farrell-Jones conjecture to hold, since for one part of the argument we use more of the HHS structure.

Asymptotically CAT(0) spaces have not been studied as much as injective spaces, but hopefully the theorem provides plenty of motivation to do so. We formulate this as a very open-ended problem:

\begin{problem}
 Study asymptotically CAT(0) spaces.
\end{problem}

\section{Further tools}
\label{sec:curtains}

In this section we discuss other tools more briefly (most of these are under development).

\subsection{Dehn-filling-like quotients}

While the theory I will describe in this subsection is still under development, it already has striking applications. In analogy with relatively hyperbolic Dehn filling \cite{Osin:perfill, GM:Dehnfill}, there are now a few known ways of taking quotients of HHGs and obtaining new HHGs. The first instance of this was in \cite{hhs_asdim}, where, in the context of mapping class groups, it is shown that quotients of mapping class groups by suitable powers of pseudo-Anosovs are still HHG. (Relatedly, random quotients will be considered in forthcoming work of Abbott, Berlyne, Ng, and Rasmussen.) In \cite{BHMS}, instead, quotients by suitable powers of all Dehn twists are shown to be HHG, and the crucial feature of this construction is that it reduces a natural measure of complexity, where minimal complexity HHG are simply hyperbolic groups. Therefore, this can be thought of as a first step in constructing hyperbolic quotients of mapping class groups, but this procedure requires certain hyperbolic groups encountered during the construction to be residually finite. Under the assumption that, say, all hyperbolic groups are residually finite, the inductive procedure is carried out in \cite{BHMS}. As it turns out, this allows to relate residual finiteness of ``enough'' hyperbolic groups to various open questions on profinite properties of mapping class groups. Taking a further, striking step in this direction, Wilton recently showed, using the construction of \cite{BHMS} almost as a black box, that if ``enough'' hyperbolic groups are residually finite then mapping class groups have the congruence subgroup property \cite{Wilton:congruence}. In a different direction, a related type of quotient of HHGs will be used in \cite{short_HHG:II} to show that many Artin groups are Hopfian. Roughly, the point for those Artin groups is that a single ``complexity reduction'' is enough to obtain a hyperbolic group. The general idea is then to exploit hyperbolic quotients to deduce algebraic properties of an HHG from those of its hyperbolic quotients.

We do not have a complete theory of ``Dehn-filling-like'' quotients of HHGs, but even partial versions of it have several potential applications, as illustrated by the examples above. I believe that a good question to drive research in this area is the following conjecture, stated in slightly simplified form for concreteness:

\begin{conj}\cite{MS:dehn_twists}
    Let $S$ be any surface of finite type, and let $g_1,\dots,g_l\in MCG(S)$. Then there exists $N\in\mathbb{N}-\{0\}$ such that $MCG(S)/\langle\langle \{g_i^{N}\}\rangle\rangle$ is hierarchically hyperbolic.
\end{conj}

The conjecture is proven for the five-holed sphere and almost all collections of elements $g_i$ in \cite{short_HHG:II}.

\subsection{Curtains}
Curtains were originally introduced for CAT(0) spaces in \cite{PSZ} as generalisations of hyperplanes in CAT(0) cube complexes. The technology has been pushed considerably further in \cite{PZ}. The construction is very general, and here we only focus on the HHS aspects of it. Roughly, a hyperbolic space can be split into two parts by fixing a geodesic, looking at the closest point projection onto it, and considering points that project before or after a certain point. In an HHS $X$ one can do something similar, first using one of the maps $\pi_Y:X\to \calC(Y)$, for $\calC(Y)$ one of the hyperbolic spaces from the HHS structure. This feeds into the very general machinery of \cite{PZ} to yield a metric space $X_\calC$ which has many remarkable properties. First, it is injective, so this gives another construction of injective metric spaces for HHSs. Secondly, it admits a bicombing by rough geodesics which interact nicely with the HHS structure, namely they are hierarchy paths. It is not clear whether the previous construction of injective spaces for HHSs had this feature. There is also another property that $X_\calC$ satisfies, namely almost convexity, a convexity property of balls. Finally, $X_\calC$ is not only a metric space, but also a median algebra. The HHS $X$ is only coarse median, so $X_\calC$ has better ``local'' properties, making it easier to work with; in some sense $X_\calC$ is at the same time the injective hull and the ``median hull'' of $X$.

Overall, besides the construction having many further applications beyond HHSs, the space $X_\calC$ sketched above is a natural and convenient geometric model to use. We note that in forthcoming work of Hagen, Petyt, and Zalloum, curtains will be used to start with HHS data (hyperbolic spaces, relations, etc.) and construct an HHS out of this. A more categorical approach to this ``inverse problem'' is being developed by Tang.

\subsection{$\mathbb R$-cubings}

It is well-known that asymptotic cones of hyperbolic groups are $\R$-trees (which can be defined as $0$-hyperbolic spaces). Actions on $\R$-trees have several applications to hyperbolic groups and beyond because of this. Casals-Ruiz, Hagen, and Kazachkov found a generalisation of $\R$-trees, which they call $\R$-cubings, and they show that asymptotic cones of HHSs are indeed $\R$-cubings, and in turn use this to show that mapping class groups, among many other hierarchically hyperbolic groups, have unique asymptotic cones up to bilipschitz equivalence \cite{R-cubings}. Roughly, an $\R$-cubing is a particular type of median space contained in a product of (typically infinitely many) $\R$-trees, and consisting of all points satisfying a certain system of equations. These equations are essentially exact versions of HHS axioms, meaning that the constants involved in those (which can be thought of as margins of error) are set to 0.

It is envisioned that one can use actions on $\R$-cubings to study certain HHGs in roughly the same way one uses $\R$-trees to study hyperbolic groups, by introducing non-discrete versions of ideas from special cube complexes in the sense of \cite{Haglund-Wise}.

\subsection{Higher-rank JSJ decompositions}

Sela is working on a series of papers studying automorphisms of colourable HHGs satisfying a weak acylindricity condition \cite{Sela1,Sela2}. This is in analogy with the theory of JSJ decompositions of hyperbolic groups, and indeed an object that can be viewed as a higher rank JSJ decomposition is constructed. The higher rank decomposition encodes the dynamics of individual automorphisms and the structure of the outer automorphism group of an HHG. For a one ended hyperbolic group $G$ there exists an epimorphism from a finite index subgroup of $Out(G)$ onto a direct product of the mapping class groups of 2-orbifolds (appearing in the canonical JSJ decomposition of $G$), where the kernel is a finitely generated virtually abelian group. The state of the art for HHGs, satisfying additional conditions, is not quite a direct analogue of this, but instead two groupoids get associated to $G$, and outer automorphisms are associated with finitely many morphisms, rather than just one. There is still a certain homomorphism from a finite index subgroup of $Out(G)$ into the direct sum of (finitely many) mapping class groups of 2-orbifolds and general linear groups. This data still suffices to show that certain natural conditions on $G$ imply that $Out(G)$ is locally virtually nilpotent, for instance, but many open questions on $Out(G)$ still remain.

The reader is also referred to \cite{BF:semisimple} for results on acylindricity conditions emerging from Sela's work.

 \section{Some open problems}

 There are many open questions on hierarchical hyperbolicity, and indeed a problem list is available at \href{http://comet.lehman.cuny.edu/behrstock/HHS.html}{http://comet.lehman.cuny.edu/behrstock/HHS.html}, together with further material on HHSs. Therefore, here I only mention the problems I advertised at my ECM talk. These arise in the context of trying to understand the topology of a manifold from a triangulation, assuming that the fundamental group is hierarchically hyperbolic. Perhaps the tools advertised in this survey could be useful to tackle some of them.

 The word and conjugacy problem are solvable for HHGs, but the proofs are rather indirect and give very slow running time, worse than exponential, in fact. It is then natural to ask:

 \begin{question}
  How quickly can one solve the word problem in an HHG?
 \end{question}

  \begin{question}
  How quickly can one solve the conjugacy problem in an HHG?
 \end{question}

 It is a very deep and difficult theorem that the isomorphism problem is solvable for hyperbolic groups \cite{Sela:iso,DG:iso}. One can then very ambitiously ask whether the same holds for all HHGs, or still ambitiously but less so, whether there are interesting sub-classes where it can be solved:

 \begin{question}
  Is the isomorphism problem solvable for HHGs? Is it solvable at least in some interesting sub-classes?
 \end{question}

 Finally, \cite{FJ} proves the Farrell-Jones conjecture for ``most'' HHGs, the main ones left out being the non-colourable examples from \cite{noncolorable}. These are however cubical groups, so they do satisfy the Farrell-Jones conjecture \cite{BL:Borel_CAT0}, and indeed I am not aware of any specific HHGs for which the Farrell-Jones conjecture is now not known to hold. Still, it is natural to ask the following, and in fact it is very likely that removing the additional assumptions from \cite{FJ} will require independently interesting mathematics.

 \begin{question}
  Does the Farrell-Jones conjecture hold for all HHGs?
 \end{question}

\bibliographystyle{alpha}
\bibliography{biblio}

\newcommand{\etalchar}[1]{$^{#1}$}
\begin{thebibliography}{GHP{\etalchar{+}}23}

\bibitem[ABD21]{ABD}
Carolyn Abbott, Jason Behrstock, and Matthew Durham.
\newblock Largest acylindrical actions and stability in hierarchically
  hyperbolic groups.
\newblock {\em Transactions of the American Mathematical Society, Series B},
  8(3):66--104, 2021.

\bibitem[ANS{\etalchar{+}}19]{ANS:UEG}
Carolyn Abbott, Thomas Ng, Davide Spriano, Radhika Gupta, and Harry Petyt.
\newblock Hierarchically hyperbolic groups and uniform exponential growth.
\newblock {\em To appear in Math. Z.}, 2019.

\bibitem[Bar24]{barak:eq_not}
Ohana Barak.
\newblock On equational noetherianity of colorable hierarchically hyperbolic
  groups.
\newblock {\em arXiv preprint arXiv:2410.00977}, 2024.

\bibitem[BF24]{BF:semisimple}
Sahana Balasubramanya and Talia Fernos.
\newblock The semi-simple theory of higher rank acylindricity.
\newblock {\em arXiv preprint arXiv:2407.04838}, 2024.

\bibitem[BHMS20]{BHMS}
Jason Behrstock, Mark Hagen, Alexandre Martin, and Alessandro Sisto.
\newblock A combinatorial take on hierarchical hyperbolicity and applications
  to quotients of mapping class groups, 2020.

\bibitem[BHS17a]{hhs_asdim}
Jason Behrstock, Mark~F. Hagen, and Alessandro Sisto.
\newblock Asymptotic dimension and small-cancellation for hierarchically
  hyperbolic spaces and groups.
\newblock {\em Proc. Lond. Math. Soc. (3)}, 114(5):890--926, 2017.

\bibitem[BHS17b]{HHS_I}
Jason Behrstock, Mark~F. Hagen, and Alessandro Sisto.
\newblock Hierarchically hyperbolic spaces, {I}: {C}urve complexes for cubical
  groups.
\newblock {\em Geom. Topol.}, 21(3):1731--1804, 2017.

\bibitem[BHS19]{HHS_II}
Jason Behrstock, Mark Hagen, and Alessandro Sisto.
\newblock Hierarchically hyperbolic spaces {II}: {C}ombination theorems and the
  distance formula.
\newblock {\em Pacific J. Math.}, 299(2):257--338, 2019.

\bibitem[BHS21]{HHS:quasiflats}
Jason Behrstock, Mark~F. Hagen, and Alessandro Sisto.
\newblock Quasiflats in hierarchically hyperbolic spaces.
\newblock {\em Duke Math. J.}, 170(5):909--996, 2021.

\bibitem[BL12]{BL:Borel_CAT0}
Arthur Bartels and Wolfgang L\"uck.
\newblock The {B}orel conjecture for hyperbolic and {${\rm CAT}(0)$}-groups.
\newblock {\em Ann. of Math. (2)}, 175(2):631--689, 2012.

\bibitem[Bon24]{bongiovanni2024extensions}
Eliot Bongiovanni.
\newblock Extensions of finitely generated veech groups, 2024.

\bibitem[Bow19]{Bowditch:quasiflats}
Brian~H Bowditch.
\newblock Quasiflats in coarse median spaces.
\newblock {\em Preprint available at homepages. warwick. ac.
  uk/masgak/papers/quasiflats. pdf}, page~48, 2019.

\bibitem[BR20]{BerlaiRobbio}
Federico Berlai and Bruno Robbio.
\newblock A refined combination theorem for hierarchically hyperbolic groups.
\newblock {\em Groups Geom. Dyn.}, 14(4):1127--1203, 2020.

\bibitem[BR22]{BR:graphs}
Daniel Berlyne and Jacob Russell.
\newblock Hierarchical hyperbolicity of graph products.
\newblock {\em Groups, Geometry, and Dynamics}, 16(2):523--580, 2022.

\bibitem[Bri10]{Bridson:MCG}
Martin~R Bridson.
\newblock Semisimple actions of mapping class groups on {CAT}(0) spaces.
\newblock {\em Geometry of Riemann surfaces}, 368:1--14, 2010.

\bibitem[Che22]{Chesser:stable}
Marissa Chesser.
\newblock Stable subgroups of the genus 2 handlebody group.
\newblock {\em Algebr. Geom. Topol.}, 22(2):919--971, 2022.

\bibitem[CRHK24]{R-cubings}
Montserrat Casals-Ruiz, Mark Hagen, and Ilya Kazachkov.
\newblock Real cubings and asymptotic cones of hierarchically hyperbolic
  groups.
\newblock {\em in preparation}, 2024.

\bibitem[DDLS24]{veech}
Spencer Dowdall, Matthew~G. Durham, Christopher~J. Leininger, and Alessandro
  Sisto.
\newblock Extensions of {V}eech groups {II}: {H}ierarchical hyperbolicity and
  quasi-isometric rigidity.
\newblock {\em Comment. Math. Helv.}, 99(1):149--228, 2024.

\bibitem[DG11]{DG:iso}
Fran{\c{c}}ois Dahmani and Vincent Guirardel.
\newblock The isomorphism problem for all hyperbolic groups.
\newblock {\em Geometric and Functional Analysis}, 21:223--300, 2011.

\bibitem[DHS17]{DHS}
Matthew~Gentry Durham, Mark~F. Hagen, and Alessandro Sisto.
\newblock Boundaries and automorphisms of hierarchically hyperbolic spaces.
\newblock {\em Geom. Topol.}, 21(6):3659--3758, 2017.

\bibitem[DHS19]{DHS:oops}
Matthew~Gentry Durham, Mark~F Hagen, and Alessandro Sisto.
\newblock Corrigendum to boundaries and automorphisms of hierarchically
  hyperbolic spaces.
\newblock {\em Geometry and Topology}, 2019.

\bibitem[DMS23]{DMS:stable}
Matthew~G Durham, Yair~N Minsky, and Alessandro Sisto.
\newblock Stable cubulations, bicombings, and barycenters.
\newblock {\em Geom. Topol.}, 27(6):2383--2478, 2023.

\bibitem[DMS25]{FJ}
Matthew~Gentry Durham, Yair Minsky, and Alessandro Sisto.
\newblock Asymptotically {CAT}(0) metrics, {Z}-structures, and the
  {F}arrell-{J}ones {C}onjecture.
\newblock {\em arXiv preprint arXiv:2504.17048}, 2025.

\bibitem[Dur23]{Durham:infinity}
Matthew~Gentry Durham.
\newblock Cubulating infinity in hierarchically hyperbolic spaces.
\newblock {\em arXiv preprint arXiv:2308.13689}, 2023.

\bibitem[DZ22]{DZ:genericity}
Matthew~Gentry Durham and Abdul Zalloum.
\newblock The geometry of genericity in mapping class groups and
  teichm$\backslash$" uller spaces via cat (0) cube complexes.
\newblock {\em To appear in Transactions of the AMS}, 2022.

\bibitem[GHP{\etalchar{+}}23]{induced_qi}
Antoine Goldsborough, Mark Hagen, Harry Petyt, Jacob Russell, and Alessandro
  Sisto.
\newblock Induced quasi-isometries of hyperbolic spaces, markov chains, and
  acylindrical hyperbolicity.
\newblock {\em Accepted in Groups, Geometry, and Dynamics}, 2023.

\bibitem[GM08]{GM:Dehnfill}
Daniel Groves and Jason~Fox Manning.
\newblock Dehn filling in relatively hyperbolic groups.
\newblock {\em Israel Journal of Mathematics}, 168:317--429, 2008.

\bibitem[Hae23]{Haettel:survey}
Thomas Haettel.
\newblock Group actions on injective spaces and helly graphs.
\newblock {\em arXiv preprint arXiv:2307.00414}, 2023.

\bibitem[Hag23]{noncolorable}
Mark~F Hagen.
\newblock Non-colorable hierarchically hyperbolic groups.
\newblock {\em International Journal of Algebra and Computation},
  33(2):337--350, 2023.

\bibitem[HHP23]{HHP:coarse}
Thomas Haettel, Nima Hoda, and Harry Petyt.
\newblock Coarse injectivity, hierarchical hyperbolicity and semihyperbolicity.
\newblock {\em Geom. Topol.}, 27(4):1587--1633, 2023.

\bibitem[HMS22]{ELTAG_HHS}
Mark Hagen, Alexandre Martin, and Alessandro Sisto.
\newblock Extra-large type artin groups are hierarchically hyperbolic.
\newblock {\em Mathematische Annalen}, pages 1--72, 2022.

\bibitem[HMS23]{converse}
Mark Hagen, Giorgio Mangioni, and Alessandro Sisto.
\newblock A combinatorial structure for many hierarchically hyperbolic spaces.
\newblock {\em arXiv preprint arXiv:2308.16335}, 2023.

\bibitem[HRSS23]{HRSS_3manifold}
Mark Hagen, Jacob Russell, Alessandro Sisto, and Davide Spriano.
\newblock Equivariant hierarchically hyperbolic structures for 3-manifold
  groups via quasimorphisms, 2023.

\bibitem[HS20]{HagenSusse}
Mark~F Hagen and Tim Susse.
\newblock On hierarchical hyperbolicity of cubical groups.
\newblock {\em Israel Journal of Mathematics}, 236(1):45--89, 2020.

\bibitem[HV24]{HV:not_automatic}
Sam Hughes and Motiejus Valiunas.
\newblock Commensurating hnn-extensions: hierarchical hyperbolicity and
  biautomaticity.
\newblock {\em Commentarii Mathematici Helvetici: A Journal of the Swiss
  Mathematical Society}, 99(2), 2024.

\bibitem[HW08]{Haglund-Wise}
Fr{\'e}d{\'e}ric Haglund and Daniel~T Wise.
\newblock Special cube complexes.
\newblock {\em Geometric and Functional Analysis}, 17:1551--1620, 2008.

\bibitem[Kar08]{Kar:thesis}
Aditi Kar.
\newblock {\em Discrete groups and {CAT}(0) asymptotic cones}.
\newblock PhD thesis, The Ohio State University, 2008.

\bibitem[Lan13]{Lang}
Urs Lang.
\newblock Injective hulls of certain discrete metric spaces and groups.
\newblock {\em Journal of Topology and Analysis}, 5(03):297--331, 2013.

\bibitem[Man24]{Mangioni:combination}
Giorgio Mangioni.
\newblock A combination theorem for hierarchically quasiconvex subgroups, and
  application to geometric subgroups of mapping class groups.
\newblock {\em arXiv preprint arXiv:2409.03602}, 2024.

\bibitem[MM99]{MM1}
Howard~A. Masur and Yair~N. Minsky.
\newblock Geometry of the complex of curves. {I}. {H}yperbolicity.
\newblock {\em Invent. Math.}, 138(1):103--149, 1999.

\bibitem[MM00]{MM2}
H.~A. Masur and Y.~N. Minsky.
\newblock Geometry of the complex of curves. {II}. {H}ierarchical structure.
\newblock {\em Geom. Funct. Anal.}, 10(4):902--974, 2000.

\bibitem[MS22]{MS:dehn_twists}
Giorgio Mangioni and Alessandro Sisto.
\newblock Rigidity of mapping class groups mod powers of twists.
\newblock {\em Accepted in Proc. R. Soc. Edinb.}, 2022.

\bibitem[MS24]{short_HHG:II}
Giorgio Mangioni and Alessandro Sisto.
\newblock Short hierarchically hyperbolic groups ii: quotients and the hopf
  property for artin groups.
\newblock {\em in preparation}, 2024.

\bibitem[Osi07]{Osin:perfill}
Denis~V Osin.
\newblock Peripheral fillings of relatively hyperbolic groups.
\newblock {\em Inventiones mathematicae}, 167(2):295--326, 2007.

\bibitem[Par24]{Park:factor_system}
Jihoon Park.
\newblock Factor system for graphs and combinatorial {HHS}.
\newblock {\em arXiv preprint arXiv:2409.08663}, 2024.

\bibitem[Pet21]{Petyt_MCG_cubical}
Harry Petyt.
\newblock Mapping class groups are quasicubical.
\newblock {\em To appear in Amer. J. Math.}, 2021.

\bibitem[PSZ24]{PSZ}
Harry Petyt, Davide Spriano, and Abdul Zalloum.
\newblock Hyperbolic models for {CAT}(0) spaces.
\newblock {\em Advances in Mathematics}, 450:109742, 2024.

\bibitem[PZ24]{PZ}
Harry Petyt and Abdul Zalloum.
\newblock Constructing metric spaces from systems of walls, with an appendix
  with {D}avide {S}priano.
\newblock {\em arXiv preprint arXiv:2404.12057}, 2024.

\bibitem[RST23]{HHS-convexity}
Jacob Russell, Davide Spriano, and Hung~Cong Tran.
\newblock Convexity in hierarchically hyperbolic spaces.
\newblock {\em Algebr. Geom. Topol.}, 23(3):1167--1248, 2023.

\bibitem[Sel95]{Sela:iso}
Zlil Sela.
\newblock The isomorphism problem for hyperbolic groups {I}.
\newblock {\em Annals of Mathematics}, pages 217--283, 1995.

\bibitem[Sel22]{Sela2}
Zlil Sela.
\newblock Automorphisms of groups and a higher rank {JSJ} decomposition {II}:
  {T}he single ended case.
\newblock {\em arXiv preprint arXiv:2209.12891}, 2022.

\bibitem[Sel23]{Sela1}
Zlil Sela.
\newblock Automorphisms of groups and a higher rank {JSJ} decomposition {I}:
  {RAAG}s and a higher rank {M}akanin-{R}azborov diagram.
\newblock {\em Geometric and Functional Analysis}, 33(3):824--874, 2023.

\bibitem[Sis19]{HHS_survey}
Alessandro Sisto.
\newblock What is a hierarchically hyperbolic space?
\newblock In {\em Beyond hyperbolicity}, volume 454 of {\em London Math. Soc.
  Lecture Note Ser.}, pages 117--148. Cambridge Univ. Press, Cambridge, 2019.

\bibitem[Vok22]{Vokes}
Kate~M Vokes.
\newblock Hierarchical hyperbolicity of graphs of multicurves.
\newblock {\em Algebraic \& Geometric Topology}, 22(1):113--151, 2022.

\bibitem[Wil24]{Wilton:congruence}
Henry Wilton.
\newblock The congruence subgroup property for mapping class groups and the
  residual finiteness of hyperbolic groups, with an appendix by {A}lessandro
  {S}isto.
\newblock {\em arXiv preprint arXiv:2410.00556}, 2024.

\bibitem[Zal23]{Zalloum_survey}
Abdul Zalloum.
\newblock Injectivity, cubical approximations and equivariant wall structures
  beyond {CAT}(0) cube complexes.
\newblock {\em arXiv preprint arXiv:2305.02951}, 2023.

\end{thebibliography}

\end{document}